\documentstyle[12pt]{article}
\begin{document}
\renewcommand{\thefootnote}{\fnsymbol{footnote}}
\newpage
\pagestyle{empty}
\setcounter{page}{0}
\renewcommand{\thesection}{\arabic{section}}
\renewcommand{\theequation}{\thesection.\arabic{equation}}
\newcommand{\sect}[1]{\setcounter{equation}{0}\section{#1}}
\vfill
\begin{center}

{\LARGE {\bf {\sf 
Modified Braid Equations for $SO_q (3)$ and noncommutative spaces
}}} \\[0.8cm]
{\large A.Chakrabarti

{\em 

Centre de Physique Th\'eorique\footnote{Laboratoire Propre 
du CNRS UPR A.0014}, Ecole Polytechnique, 91128 Palaiseau Cedex, France.\\
e-mail chakra@cpht.polytechnique.fr}}

\end{center}

\smallskip

\smallskip 

\smallskip

\smallskip

\smallskip

\smallskip 

\begin{abstract}
General solutions of the $\hat{R}TT$ equation with a maximal number of free
parameters in the specrtal decomposition of vector $SO_q (3)$ $\hat{R}$
matrices are implemented to construct modified braid equations $(MBE)$.
These matrices conserve the given, standard, group relations of the nine
elements of $T$, but are not constrained to satisfy the standard braid
equation ($BE$). Apart from $q$ and a normalisation factor our $\hat{R}$
contains two free parameters, instead of only one such parameter for
deformed unitary algebras studied in a previous paper $[1]$ where the
nonzero right hand side of the $MBE$ had a linear term proprotional to
$(\hat{R}_{(12)} - \hat{R}_{(23)})$. In the present case the r.h.s. is, in
general, nonliear. Several particular solutions are given (Sec.2) and the
general structure is analysed (App.A). Our formulation of the
problem in terms of projectors yield also two new solutions of standard
(nonmodified ) braid equation (Sec.2) which are further discussed 
(App.B). The noncommutative $3$-spaces obtained by implementing such
generalized $\hat{R}$ matrices are studied (Sec.3). The role of coboundary
$\hat{R}$ matrices ( not satisfying the standard $BE$ ) is explored. The
$MBE$ and Baxterization are presented as complementary facets of the same
basic construction, namely, the general solution of $\hat{R}TT$ equation
(Sec.4). A  new solution is presented in this context. As a simple but
remarkable particular case a nontrivial solution of $BE$ is obtained
(App.B) for $q=1$. This solution has no free parameter and is not
obtainable by twisting the identity matrix. In the concluding remarks
(Sec.5), among other points, generalisation of our results to $SO_{q}(N)$
is discussed.   
\end{abstract}

\vfill
\noindent
CPHT-S 011.0301
\noindent
\newpage 

\pagestyle{plain}

\section{Introduction}

  In a previous paper $[1]$ a particular class of inhomogeneous, modified,
braid equations $(MBE)$ was shown to correspond to general solutions of
$\hat{R}TT$ relations. Fundamental $2\times 2$  $T$ matrices and the
corresponding $4\times 4$  $\hat{R} (= PR)$ matrices for $GL_{p,q}(2)$,
$GL_{g,h}(2)$ and $GL_{q,h}(1/1)$ were used as examples. The inputs were
the known group relations of the elements $(a,b,c,d)$ for each of the above 
cases. Then the most general solution ( without imposing the Yang Baxter
equation for $R$ or, equivalently, the braid equation for $\hat{R}$) was
sought, for each case, of the relation   
     
\begin{equation}
          \hat{R}T_{1} T_{2} = T_{1} T_{2}\hat{R}
\end{equation}

where    
         $$ T_{1} = T \otimes I_{2},\quad  T_{2} = I_{2} \otimes T$$ 

 The only constraint on $\hat{R}$ was the conservation of the given group
relations for $(a,b,c,d)$. In each case $\hat{R}$ was found to depend, ,
apart from the two parameters $\bigl( (p,q),(g,h),(q,h) \bigr)$ linearly
on a third one $(K)$ such that, for a suitable normalisation, one obtains 

\begin{equation}
            \hat{R}_{(12)}\hat{R}_{(23)}\hat{R}_{(12)} -
\hat{R}_{(23)}\hat{R}_{(12)}\hat{R}_{(23)} = \biggl( { K\over {
K_{1}}} - 1
\biggr) \biggl( {K\over { K_{2}}} - 1 \biggr) \biggl(\hat{R}_{(23)} -
\hat{R}_{(12)}\biggr) 
\end{equation}

This is the $MBE$ with  
\begin{equation}
(K_{1},K_{2}) = (1,p/q),\quad  (1,1),\quad (1,1/q)
\end{equation}
respectively for the above-mentioned cases.

 It was pointed out in $[1]$ that $(1.2)$ reexpressed in terms of $R$,
provides a particular, interesting class of modified quantum $YB$ equations
$(MQYBE)$ of Gerstenhaber,Giaquinto and Schack ( see [2] and sources cited
therein ).

      The special features of $(1.2)$, as indicated in $[1]$ are as follows.

$(1)$: The explicit stucture on the right carries interesting
information. After obtaining the general solution of $(1.1)$ and the $MBE$
it corresponds to one obtains the unmodified,standard braid or $YB$
matrices as byproducts. One just sets $K=K_{1}$ or $K=K_{2}$, the two
solutions being related through

\begin{equation}
   \hat{R}(K_{2}) =  {\hat{R}(K_{1})}^{-1}
\end{equation}

$(2)$: Setting, as explained in $[1]$, 

\begin{equation}
   K = 2K_{1}K_{2}{(K_{1} +K_{2})}^{-1}
\end{equation}
one obtains

\begin{equation}
   {(\hat{R}(K))}^{2} =  I
\end{equation}

Hence the construction of "triangular" ( or "unitary" or "coboundary") $R$
matrices is again reduced to the choice of a particular value of $K$ in
the general solution. This aspect will be studied further below.

$(3)$: It was pointed out in $[1]$ that  $MBE$ and Baxterization are
complementary facets of the same basic construction , namely, the
general solution of $(1.1)$. This links $MBE$ to integrable models. This
aspect will be taken up in $Sec.4$ with new examples.

$(4)$: The parameter $K$ plays an interesting role in noncommutative spaces
obtained by implementing the general $\hat{R}(K)$. The detailed discussion
in $[1]$ of this aspect will be generalized to higher dimensions in
$Sec.(3)$.
\smallskip

\smallskip

   Spectral decomposition and generalization to higher dimensions:

\smallskip

\smallskip

Our construction can be generalized to higher dimensional cases most
conveniently by introducing arbitrary constant coefficients in the spectral
decomposition of $R$ matrices for {\it {vector representations}}. The
standard solutions ( not "modified" in our sense ) are well-known $[3,4]$.
Instead of wading through larger and large number of group relations (
$n^{2}(n^{2} - 1)/2$ for $n^2$ elements of $T$) one starts with the
following results for vector representations.

        For $GL_{q}(N)$ one has, in terms of the projectors $P^{(\pm)}$, 

\begin{equation}
\hat{R}=qP^{(+)}-q^{-1}P^{(-)}
\end{equation}
where
\begin{equation}
P^{(+)}P^{(-)}=0, \qquad \bigl( P^{(\pm)}\bigr)^{2}=P^{(\pm)}, \qquad
P^{(+)}+P^{(-)}=I
\end{equation}

Here $\hat{R}$ satisfies the braid equation and
\begin{equation}
\bigl(\hat{R}-qI\bigr)\bigl(\hat{R}+q^{-1}I\bigr)=0
\end{equation}

One has also
\begin{equation}
P^{(+)}=\frac{\bigl(\hat{R}+q^{-1}I\bigr)}{\bigl(q+q^{-1}\bigr)},\qquad
P^{(-)}=-\frac{\bigl(\hat{R}-qI\bigr)}{\bigl(q+q^{-1}\bigr)}
\end{equation}

If one sets, {\it {with the same projectors}},

\begin{equation}
\hat{R}(u,v)=uP^{(+)}+vP^{(-)}
\end{equation}
where $(u,v)$ are non-zero, unequal but otherwise arbitrary parameters, one
obtains
\begin{equation}
\bigl(\hat{R}(u,v)-uI\bigr)\bigl(\hat{R}(u,v)-vI\bigr)=0
\end{equation}

  \begin{equation}
P^{(+)}=\frac{\bigl(\hat{R}(u,v)-vI\bigr)}{\bigl(u-v\bigr)},\qquad
P^{(-)}=\frac{\bigl(\hat{R}(u,v)-uI\bigr)}{\bigl(v-u\bigr)}
\end{equation}

Of the two parameters $(u,v)$ one can be fixed by choosing a suitable
normalisation, leading effectively to one independent, arbitrary parameter.
Apart form differences in notations our construction of $\hat{R}(K)$ in
$[1]$ ( see Sec.3.2 in $[1]$ ), namely
\begin{equation}
\hat{R}(K)=(1-(K/K_{1} + K/K_{2}))P_{1}+P_{2}
\end{equation}
corresponds directly to $(1.11)$ above. This and certain other aspects of
our previous formalism {\it {can be directly generalised }}to $GL_{q}(N)$,
though the noncommutative spaces will now be of $N$ dimensions.

    For $SO_{q}(N)$ ( and for $Sp_{q}(N)$ which we do not consider here )
there is a major change. One has now {\it {three}} projectors in the
spectral decomposition of $\hat{R}$ matrices for vector representations.
The consequences for $MBE$ will be seen to be important.

     For $\hat{R}$ satisfying the braid equation one obtains $([3],[4])$

\begin{equation}
\bigl(\hat{R}-qI\bigr)\bigl(\hat{R}+q^{-1}I\bigr)\bigl(\hat{R}-q^{1-N}I\bigr)
=0
\end{equation}
and 

\begin{equation}
\hat{R}=qP^{(+)}-q^{-1}P^{(-)}+q^{1-N}P^{(0)}
\end{equation}

where ( with $(i,j)$ denoting a pair from $(+,-,0)$ )
\begin{equation}
P^{(i)}P^{(j)}=P^{(i)}{\delta}_{ij}, \qquad P^{(+)}+P^{(-)}+P^{(0)}=I
\end{equation}

Generalising as before we introduce ( with non-zero and unequal $(u,v,w)$
and the projectors being the {\it  same} as before, independent of $(u,v,w)$
)

\begin{equation}
\hat{R}(u,v,w)=uP^{(+)}+vP^{(-)}+wP^{(0)}
\end{equation}

Now ( denoting $\hat{R}(u,v,w)$ as $\hat{R}$ ),
\begin{equation}
\bigl(\hat{R}-uI\bigr)\bigl(\hat{R}-vI\bigr)\bigl(\hat{R}-wI\bigr)=0
\end{equation}
and
$$ 
P^{(+)}=\frac{\bigl(\hat{R}-vI\bigr)\bigl(\hat{R}-wI\bigr)}
{\bigl(u-v\bigr)\bigl(u-w\bigr)}
$$
\begin{equation}  
P^{(-)}=\frac{\bigl(\hat{R}-uI\bigr)\bigl(\hat{R}-wI\bigr)}
{\bigl(v-u\bigr)\bigl(v-w\bigr)}
\end{equation}
$$
P^{(0)}=\frac{\bigl(\hat{R}-uI\bigr)\bigl(\hat{R}-vI\bigr)}
{\bigl(w-u\bigr)\bigl(w-v\bigr)}
$$

Here, fixing the normalisation, {\it two} independent parameters are left.

 In the next section we will study the $MBE$ corresponding to $(1.18)$. We
will set $N=3$. This will permit us to display explicitly matrices of
managaeble size. {\it {The essential new features will, however, be present
already for}} $N=3$. 

Let us now note how the number of coboundary (or unitary) solutions for
vector representations changes with the number of projectors.

    For $GL_{q}(N)$ it is seen from $(1.8)$ and $(1.11)$ that
\begin{equation}
 \Bigl({\hat{R}(u,v)}\Bigr)^{2}=u^{2}P^{(+)}+v^{2}P^{(-)} =P^{(+)}+P^{(-)}
\end{equation}
for
$$ u^2 = v^2 = 1.$$ 

 Hence, apart from an overall $(\pm)$ sign, the only nontrivial solution is

\begin{equation}
\hat{R}_{c} = P^{(+)}-P^{(-)}= I-2P^{(-)}=-(I-2P^{(+)}) ;\qquad
{\hat{R}_{c}}^{2}=I
\end{equation}

For $SO_{q}(n)$, from $(1.17)$ and $(1.18)$ apart from an overall sign one
obtains analogously {\it{three}} solutions
\begin{equation}
\hat{R}_{c} = (I-2P^{(+)}); \quad(I-2P^{(-)}); \quad  (I-2P^{(0)})
\end{equation}

Each satisfies 
$$ {\hat{R}_{c}}^{2}=I$$

and the product of any two gives the third one with a change of sign. Thus,
for example
\begin{equation}
 (I-2P^{(+)})(I-2P^{(-)})=- (I-2P^{(0)})
\end{equation}

 If the coefficient $(-2)$ in $(1.23)$ is replaced by an arbitrary number
$\hat{R}$ still satifies as is easily seen a {\it{quadratic}} equation, not
the cubic $(1.19)$.

  If complex solutions are considered for real $q$ but with complex
coefficients in $(1.18)$ , one obtains the unitarity relation when $(u,v,w)$
are phases  in $(1.18)$. Thus ( with real deltas ) 
\begin{equation}
\hat{R}=e^{i\delta_{1}}P^{(+)}+e^{i\delta_{2}}P^{(-)}+e^{i\delta_{3}}P^{(0)}
\end{equation}
gives, since the projectors are symmetric for the orhogonal case, 
$${\hat{R}}^{\dagger}\hat{R}= I$$

  In $(1.18)$, $(u,v,w)$ were postulated to be unequal. This permits one
to express {\it {all the three}} projectors in terms of $\hat{R}$ as in
$(1.20)$. But this is not obligatory. As noted in $(1.23)$ other cases (
with $u=v =-w=1$ and so on ) can indeed be of special interest. For
$(1.23)$ in each case one has {\it{two}} mutually orthogonal combinations.
Selecting the second case, for example, one obtains
\begin{equation}
P^{(-)}= -{\frac{1}{2}}(\hat{R}_{c}-I), \qquad  P^{(+)}+P^{(0)} =
{\frac{1}{2}}(\hat{R}_{c}+I)
\end{equation}

 We conclude with a fully explicit statement of our approach. If one has

$$ \hat{R}T_{1}T_{2}= T_{1}T_{2}\hat{R}$$
then any function $f(\hat{R})$ of $\hat{R}$ satisfies
\begin{equation}
f(\hat{R})T_{1}T_{2}= T_{1}T_{2}f(\hat{R})
\end{equation}

 To start with let us suppose ( for definiteness,  such a starting
point not being essential ) that $\hat{R}$ satifies $BE$.

      For $GL_{q}(N)$, $\hat{R}$ satisfies a quadratic constraint $(1.9)$.
Hence any power series in $\hat{R}$ can be reduced to a linear function in 
 $\hat{R}$. Hence, apart from an overall normalisation factor, the most
general solution of $(1.27)$, for a given set of group relations, becomes
effectively $(1.14)$ as studied in $[1]$. ( It is easy to see from $(1.7)$
and $(1.10)$ that even fractional powers of $\hat{R}$ can be obtained as a
linear funtion of $\hat{R}$ but, in general, with complex coefficients.
Having noted this, we will usually implicitly consider real coefficients. An
analogous situation will hold for the orthogonal case considered below.
Complex coefficients, such as in $(1.25)$ will not be introduced
explicitly. Except when roots of unity are involved, complexification of our
formalism is however straightforward.)

 For $SO_{q}(N)$, $\hat{R}$ satisfies the a cubic constraint $(1.15)$.
Hence the general solution, using analogous arguments, is seen effectively
to be, with constant coefficients $c_{i}$,
\begin{equation}
f(\hat{R})=c_{1}(\hat{R})^{2} + c_{2}\hat{R}+c_{3}I
\end{equation}

Using the spectral decomposition $(1.16)$ along with $(1.17)$,
$$
f(\hat{R}) = c_{1}(q^{2}P^{(+)}+q^{-2}P^{(-)}+q^{2(1-N)}P{(0)})
$$
$$
+c_{2}(qP^{(+)}-q^{-1}P^{(-)}+q^{(1-N)}P{(0)})
$$
$$
+c_{3}(P^{(+)}+P^{(-)}+P^{(0)})
$$

Hence, collecting together the coefficients one obtains the form

\begin{equation}
f(\hat{R})= uP^{(+)}+vP^{(-)}+wP^{(0)}
\end{equation}

This is the motivation for $(1.18)$. The starting point is the most
general solution, for a given set of group relations, as given by
$(1.29)$. The right hand side of the $MBE$ will be a {\it{consequence}} (
see App.A). In this larger space one then looks for points with particularly
attractive properties ( for example, those corresponding to $\hat{R}_{c}$ )
 and, more generally, explores the consequences of the free parameters in
$(1.29)$ such as in related noncommutative geometries (Sec.3).

\section{MBE for $SO_{q}(3)$:}

 We fix the normalisation by choosing the top left element (row $1$,
col.$1$) to be unity. In order to simplify the explicit form of $\hat{R}$
we denote the remaining parameters in $(1.18)$ as follows  

\begin{equation}
\hat{R}= P^{(+)}+(1+a(1+q^{2}))P^{(-)}+(1+b(1+q+q^{2}))P^{(0)}
\end{equation}
$$
=I+a(1+q^{2})P^{(-)}+b(1+q+q^{2})P^{(0)}
$$
The projectors are given explicitly at the end of $App.A$. The $9\times 9$
symmetric $\hat{R}$ is now
$$
\hat{R}(a,b;q)= 
$$
\begin{equation}{\scriptsize
\pmatrix{
1 &0 &0 &0 &0 &0 &0 &0 &0 \cr
0 &(1+a) &0 &-aq &0 &0 &0 &0 &0 \cr
0 &0 &(1+aq+b) &0 &(b+a(q-1))\sqrt{q} &0 &(b-a)q &0 &0 \cr
0 &-aq &0 &(1+aq^{2}) &0 &0 &0 &0 &0 \cr
0 &0  &(b+a(q-1))\sqrt{q} &0 &(1+a(q-1)^{2}+bq) &0 &(bq-a(q-1))\sqrt{q} &0
 &0 \cr
0 &0 &0 &0 &0  &(1+a) &0 &-aq &0 \cr
0 &0 &(b-a)q &0 &(bq-a(q-1))\sqrt{q} &0 &(1+aq+bq^{2}) &0 &0 \cr
0 &0 &0 &0 &0 &-aq &0 &(1+aq^{2}) &0  \cr
0 &0 &0 &0 &0 &0 &0 &0 &1                                                  
 }.}
\end{equation}

This $\hat{R}(a,b;q)$  satisfies the braid equation for 

$$ (1):a=-q^{-2},\quad  b= -q^{-2}+q^{-3}$$
\begin{equation}
(2): a =-1,\qquad  b=-1+q
\end{equation}

the two sets giving mutually inverse matrices. 

For this $\hat{R}$ ( the parameters $(a,b)$ being implicit and $I$ being the
$3\times 3$ unit matrix ) we define

\begin{equation}
\hat{R}_{12}=\hat{R}\otimes{I}, \quad \hat{R}_{23}={I}\otimes\hat{R}
\end{equation}
 
The general srtucture of the $MBE$ is presented in $App.A$. We present here
three cases obtained for particular constraints on $(a,b)$. One has

\begin{equation}
\hat{R}_{12}\hat{R}_{23}\hat{R}_{12}
-\hat{R}_{23}\hat{R}_{12}\hat{R}_{23} = l_{1}(\hat{R}_{12}- \hat{R}_{23})
+l_{2}({\hat{R}_{12}}^{2}- {\hat{R}_{23}}^{2}) 
\end{equation}

where for

$Case 1$:  $a=0$ and arbitrary $b$,

\begin{equation}
l_{2} =0; \qquad l_{1}= (1+b(1+q+q^{2})+b^{2}q^{2})
\end{equation}

$Case 2$:  $b=(1-q)a$

\begin{equation}
l_{2} =(1+a); \qquad l_{1}= (3+2(2+q^{2})a+(1+2q^{2})a^{2})
\end{equation}

Setting $a=-1$ one obtains the case $(2)$ of $(2.32)$ with 
$$l_{1}=0,\quad l_{2}=0$$

$Case 3$:  $b=(1-q^{-1})a$

\begin{equation}
l_{2} =(1+q^{2}a); \qquad l_{1}= (3+2(1+2q^{2})a+q^{2}(2+q^{2})a^{2})
\end{equation}

Setting $a=-q^{-2}$ one obtains the case $(1)$ of $(2.32)$ with 
$$l_{1}=0,\quad l_{2}=0$$

Let us note the following features:

 Thr right hand side of $(2.34)$ is linear only for $a=0$ (Case $1$). This
is evidently not included in the standard cases $(2.32)$.{\it {Yet the
braid equation is satisfied}} for

$$ a=0,\qquad  (b^{2}q^{2}+b(1+q+q^{2}) +1)=0$$

hence ( when $q \neq 0$ ) for 

\begin{equation}
b= -\frac {1}{2q^{2}}(1+q+q^{2}) \pm \frac {1}{2q^{2}}{\bigl ((1+3q+q^{2})
(1-q+q^{2})\bigr )}^{\frac{1}{2}}
\end{equation}

This gives real $b$ for $q>0$. (See $App.B$ for further discussion.)

To complete the picture,we note that the braid matrix becomes for
$$q=0, \quad a=0, \quad b=-1$$
\begin{equation}
\hat{R}(0,-1;0) = diag (1,1,0,1,1,1,1,1,1 )
\end{equation}

When $a$ and $b$ are independent and arbitrary even the quadratic terms of
the r.h.s. of $(2.34)$ do not suffice ($App.A$).

\section{Noncommutative $3$-space from $\hat{R}(a,b;q)$:}

 When $\hat{R}$ satisfies the braid equation (for $(2.32)$) the quantum
vector space is discussed in $[3,4,5]$. (See in particular $Sec.(9.3.2) $
of $[4]$ and $Ex.(4.1.22)$ of $[5]$. Our results below are to be compared
to these treatments.) We will treat the more general case with parameters
$(a,b)$. The explicit form of the r.h.s. of the $MBE$ ($Sec.2$) is not
directly relevant here. We treat $(a,b)$ as free parameters to start
with. With a slight change of notation ( with respect to $Sec.1$ ) we set in
$(2.30)$  

$$v=a(1+q^{2}); \quad w=b(1+q+q^{2})$$
giving

\begin{equation}
\hat{R}= I+vP^{(-)}+wP^{(0)} \quad = P^{+}+(1+v)P^{(-)}+(1+w)P^{(0)} 
\end{equation}

The coordinates are denoted  $(x_{-},x_{0},x_{+})$. Let $(x\otimes {x})$
( without "tilde" for simplicity) denote the $9$-component {\it {column}}
obtained from the tensor product. Let $({\xi}_{-},{\xi}_{0},{\xi}_{+})$
denote the differentials $(dx_{-},dx_{0},dx_{+})$. Let the columns for the
other tensor products be denoted, in evident notations, as $({\xi}\otimes
{\xi})$,  $(x\otimes {\xi})$,  $({\xi}\otimes {x})$.

As in $[1]$ we will adopt prescriptions that give commutators of
$(x_{i},x_{j})$ and of $({\xi}_{i},{\xi}_{j})$ independent of $(v,w)$
while those of $(x_{i},{\xi}_{j})$ do depend on them.

Let 
$$(\hat{R}-I)(\hat{R}-(1+w)I)(x\otimes {x})=0$$
or
 \begin{equation}
P^{(-)}(x\otimes {x})=0 
\end{equation}

This agrees with $[3,4,5]$. Now set

 \begin{equation}
(x\otimes {\xi})=M({\xi}\otimes {x}) 
\end{equation}

Exterior derivation gives

 \begin{equation}
({\xi}\otimes {\xi})=-M({\xi}\otimes {\xi}) 
\end{equation}

or $$ (M+I)({\xi}\otimes {\xi})=0$$.

Now exterior derivation of $(3.41)$ along with $(3.42)$ gives the typical
constraint

\begin{equation}
(\hat{R}-I)(\hat{R}-(1+w)I)(M+I)=0
\end{equation}

Hence one can choose ($k$ being an arbitrary constant parameter)
\begin{equation}
(M+I)=k(\hat{R}-(1+v)I)=k(-vP^{(+)}+(w-v)P^{(0)})
\end{equation}

From $(3.43)$ and $(3.45)$, due to the orthogonality of the projectors,
one obtains ( in agreement with $[4]$ and $[5]$ )
\begin{equation}
P^{(+)}({\xi}\otimes {\xi})=0,\quad P^{(0)}({\xi}\otimes {\xi})=0
\end{equation}

From $(3.41)$ and $(3.46)$ one obtains ( as in the standard treatments cited
above ) 
$$x_{-}x_{0}=qx_{0}x_{-}, \qquad x_{0}x_{+}=qx_{+}x_{0}$$
\begin{equation}
x_{+}x_{-}-x_{-}x_{+}=h{x_{0}}^{2} 
\end{equation}
with  $$h\equiv \Bigl(\sqrt{q}-{\frac{1}{\sqrt{q}}}\Bigr) $$
and

$${{\xi}_{-}}^{2}=0,\quad {{\xi}_{+}}^{2}=0,\quad
{\xi}_{-}{\xi}_{+}+{\xi}_{+}{\xi}_{-}=0 $$

$$q{\xi}_{-}{\xi}_{0}+{\xi}_{0}{\xi}_{-}=0,\quad
q{\xi}_{0}{\xi}_{+}+{\xi}_{+}{\xi}_{0}=0$$

\begin{equation}
 {{\xi}_{0}}^{2}=h{\xi}_{-}{\xi}_{+}
\end{equation}

Now we come to the part specific to our formalism. We define

$${\Phi}_{-}=({\xi}_{-}x_{0}-q{\xi}_{0}x_{-}), \quad
{\Phi}_{+}=({\xi}_{0}x_{+}-q{\xi}_{+}x_{0})$$
\begin{equation}
 {\Phi}_{0}=({\xi}_{-}x_{+}+{\sqrt{q}}{\xi}_{0}x_{0}+q{\xi}_{+}x_{-}), \quad
{\Phi}_{0}'=({\xi}_{-}x_{+}+h{\xi}_{0}x_{0}-{\xi}_{+}x_{-})
\end{equation}

Then, implementing the definition of $M$ in terms of $\hat{R}(a,b)$ and the
explicit form of the latter one obtains from $(3.42)$, denoting
$$k_{1}=-(kv+1)=-(ka(q^{2}+1)+1),$$
the module srtucture
$$x_{-}{\xi}_{-}=k_{1}{\xi}_{-}x_{-}; \quad
x_{-}{\xi}_{0}=k_{1}{{\xi}_{-}}x_{0}+ka{{\Phi}_{-}}$$ 
$$x_{-}{\xi}_{+}=k_{1}{{\xi}_{-}}x_{+}+kaq{{\Phi}_{0}'}+kb{\Phi}_{0}; \quad
x_{0}{\xi}_{-}=k_{1}{{\xi}_{0}}x_{-}-kaq{{\Phi}_{-}} $$
$$x_{0}{\xi}_{0}=k_{1}{{\xi}_{0}}x_{0}+ka(q-1){\sqrt q}
{{\Phi}_{0}'}+kb{\sqrt q}{\Phi}_{0}$$ 
$$x_{0}{\xi}_{+}=k_{1}{{\xi}_{0}}x_{+}+ka{{\Phi}_{+}}; \quad
x_{+}{\xi}_{-}=k_{1}{{\xi}_{+}}x_{-}-kaq{{\Phi}_{0}'}+kbq{\Phi}_{0}$$
\begin{equation}
x_{+}{\xi}_{0}=k_{1}{{\xi}_{+}}x_{0}-kaq{{\Phi}_{+}};\quad
x_{+}{\xi}_{+}=k_{1}{\xi}_{+}x_{+}
\end{equation}

One can verify that one obtains the relations given in $Ex.4.1.22$  of
$[5]$ ( page $133$) on setting
\begin{equation}
k=q^{2},\quad a=-q^{-2},\quad b=-q^{-2}+q^{-3}
\end{equation}

For 
\begin{equation}
k=1,\quad a=-2{(1+q^{2})}^{-1},\quad b=0
\end{equation}

one obtains the case $(1.26)$ of $\hat{R}_{(c)}$ ( with
${\hat{R}_{(c)}}^{2}=I $) where
\begin{equation}
P^{(-)}={\frac{1}{2}}(\hat{R}_{(c)}-1); \quad (M+I)=2(P^{(+)}+P^{(0)})
\end{equation}

Hence $(3.47)$ and $(3.48)$ are conserved along with a particularly simple
form of $(3.50)$. Here one moves out of the restricted space of solutions of
$BE$ ( or $YBE)$ ) to implement the particular simpicity of
$\hat{R}_{(c)}$. (See the relevant remarks in $Sec.5)$.)

\section{MBE and Baxterization:}

  In $[1]$ we briefly pointed out that $MBE$ and Baxterization are two
complementary aspects of the same basic construction: the general solution
of $\hat{R}TT$ equation for a given set of group relations of the
elements of $T$. For the cases considered in $[1]$ ( generalisable to
$GL_{q}(N)$ ) the correspondance is relatively simple. In $(1.2)$ the same,
single parameter $K$ appears in each member on the left leading to the
non-zero r.h.s. ( thus modifying the $BE$ ) as shown in $(1.2)$. In a
complementary approach, one can vary $K$ in different members on the left in
a prescribed fashion (indicated in $[1]$ ) so that the r.h.s. remains zero.
This is Baxterization. The {\it same} parameter that leads to $MBE$ thus
leads also to integrable systems in a complementary fashion.

  One can make a parallel study for $SO_{q}(3)$ ( generalisable to
$SO_{q}(N)$ and $Sp_{q}(N)$ ). But the presence of three projectors and
hence ( apart from a normalisation factor )  of two arbitrary parameters
leads to a more complex situation. Even restricted cases give the $MBE$ of
$(2.34)$ with quadratic terms on the right, the general structure being
given in $App.A$. Let us now look at the complementary situation, namely,
Baxterization.

 In $[4]$ a solution is given ( p.$295-297$ ) of
\begin{equation}
\hat{R}_{12}(x)\hat{R}_{23}(z)\hat{R}_{12}(y) -
\hat{R}_{23}(y)\hat{R}_{12}(z)\hat{R}_{23}(x) = 0
\end{equation}

for the restricted case where

\begin{equation}
z=xy
\end{equation}

In our notations the solution of $[4]$ reads

\begin{equation}
\hat{R}(x)= I +
\frac{x-1}{q-q^{-1}}\bigl(qI-(q+q^{-1})P^{(-)}-(q-q^{-2})P^{(0)}\bigr) +
{\frac{x+1}{x{\alpha}_{(\pm)} +1}}
\biggl(1-{\frac{q^{4}-q^{-4}}{q-q^{-1}}}\biggr)P^{(0)}
\end{equation}
where for $SO_{q}(3) \quad  {\alpha}_{(\pm)}=\pm q^{(2\pm 1)}$.

 ( We have suppressed an overall factor $h(x)$ which cancels out in
$(4.54)$. We have used $(81)$ of $p.275$ of $[4]$ to express $K$ of this
reference by $P^{(0)}$. Finally we have written $x$ for $x^{\gamma}$ of
$[4]$. One can rewrite  $(4.55)$ as $z^{\gamma}=x^{\gamma}y^{\gamma}$
and then redefine again absorbing $\gamma$. Our notation displays the
single parameter that is effectively implemented. Adjusting the
normalisation one finds a particular case of $(3.40)$.)

  The restriction $(4.55)$ is however not essential in Baxter's criterion
for commuting tansfer matrices. One can ask whether for given $(x,y)$ a
$z$ can be found assuring $(4.54)$. We present a relatively simple
solution providing the complementary facet of the $MBE$ of $(2.35)$. Let 

\begin{equation}
\hat{R}(w)=I+wP^{(0)}
\end{equation}

Then 

\begin{equation}
\hat{R}_{12}(w)\hat{R}_{23}(w')\hat{R}_{12}(w'') -
\hat{R}_{23}(w'')\hat{R}_{12}(w')\hat{R}_{23}(w) = 0
\end{equation}

for 
\begin{equation}
w'= \frac{w+w''+ww''}{1-q^{2}{(1+q+q^{2})}^{-2}ww''}
\end{equation} 

We have used results in $App.A$ ( in particular $(6.74)$ ) to derive
$(4.59)$, but it can be verified directly. Generalisations are not evident.
But a systematic study of possibilities in this context is desirable.    

\section{Remarks:}

    We conclude by noting the following points.

$(1):$  After the introductory remarks $(Sec.1)$ on the spectral
decomposition of $\hat{R}$ for vector representations of $GL_{q}(N)$,
$SO_{q}(N)$ and $Sp_{q}(N)$, from $Sec.2$ onwards we restricted our study
to $SO_{q}(3)$. But a substantial part of our results are evidently
generalisable to $SO_{q}(N)$ with $N>3$. The crucial feature is the number
of projectors in the spectral decomposition. The $MBE$ for  $GL_{q}(2)$ with
two projectors and that for $SO_{q}(3)$ with three  exhibit major
differences, made explicit here. But in $SO_{q}(N)$ the number of
projectors does not vary with $N$. Still a careful stdy of of the case of
$SO_{q}(4)$ and comparison of the results with those for $SO_{q}(3)$ would
be of real interest. This is beyond the scope of this paper. 

$(2):$ In $[1]$ we started with the criterion of using the most general
solution $\hat{R}$ of the $\hat{R}TT$ relations for a given set of group
relations of the elements of $T$. This, being implemented in the standard
trilinear structure of the braid equation, modified the right hand making
it non-zero but {\it {linear}} in the $R$'s as shown in $(1.2)$. {\it{No
apriori postulate was was made concerning the r.h.s. of the equation.}} The
explicit form was a {\it{consequence}} of the free parameter $K$ in
$\hat{R}$. It was then noted in $[1]$ that the
$MBE$ thus obtained ( eqn.$(1.2)$ of this paper ) coincided with that
introduced in $[2]$  ( and sources cited there ). Now this is seen to be a
{\it {coincidence}} valid for the cases studied in $[1]$ ( generalisable
to $GL_{q}(N)$ ). All those cases involved two projectors ( the sum being
$I$ ). As soon as this number increases ( such as already for $SO_{q}(3)$
) the r.h.s. has a more complex structure. Our starting point ( the
general solution of $\hat{R}TT$ ) is exactly the same here. But only in the
very particular case $(2.35)$ one has a linear srtucture on the right.

$(3):$ In $Sec.3$ we present all the relations involving the coordinates
and the differentials. But much remains to be done to better understand the
noncommutative space thus obtained. The properties of the $\Phi$'s
introduced in $(3.49)$ deserve  study. To render the geometry
more transparent one should construct the "frame basis" in terms of
operators commuting with the algebra. ( See $[5],[6]$ and sources cited
there.) Thus equipped, one can study possible attractive consequences
concerning the metric of implementing $\hat{R}_{c}$ as in $(3.53)$. We
hope to explore these aspects elsewhere.

$(4):$  After $[1]$ we emphasise here again the complementary nature
of $MBE$ and Baxterization. We provide a simple new example of the
latter in $(4.57)$. Let us repeat another point made in $[1]$. The
standard braid equation can be made to correspond to the third
Reidemeister move in knot theory. Hence $(4.54)$ can be viewed as a
parametrisation of this move. It would be worth exploring in the
context of knot theory whether this provides access to a richer class
of invariants.

$(5):$ Here we have studied coboundary $\hat{R}$ martices in vector
representations in terms of projectors. An approach using Drinfeld's
transformation can be found in $[7]$.

$(6):$ Finally let us recapitulate the remarkable properties of the
solution provided by $(2.38)$. It satifies the $\hat{R}TT$ equation for the
standard group relations of $SO_{q}(3)$. It satisfies the standard  braid
equation ($(2.34)$ with zero r.h.s.). It continues ( as shown in $App.B$) to
be nontrivial even for $q=1$. It can be relatively simply Baxterized as
shown at the end of $Sec.4$. This sets the stage for a full study of the
corresponding integrable model. This also suggests a search of analogous new
solutions, more generally, for $SO_{q}(N)$ and $Sp_{q}(N)$ and also for
higher dimensional representations of $SO_{q}(3)$. We hope to explore such
possibilities elsewhere. 

\section {APPENDIX A : General structure of MBE}

 In the notation of $(2.30)$ or $(3.40)$

\begin{equation}
\hat{R}= I+vP^{(-)}+wP^{(0)} 
\end{equation}
For elucidating the sructure of the consequent $MBE$  we start with a
number of definitions and auxiliary relations .
 
 We define

\begin{equation} 
 X_{1}=P^{(-)}_{(12)},\quad X_{2}=P^{(-)}_{(23)};  \quad
Y_{1}=P^{(0)}_{(12)}, \quad  Y_{2}=P^{(0)}_{(23)}
\end{equation} 

 The orthonormal properties of the projectors imply ( for $i=(1,2)$ )

\begin{equation}
{X_{i}}^{2}= X_{i},\quad {Y_{i}}^{2}= Y_{i}; \qquad
X_{i}Y_{i}= Y_{i}X_{i}=0
\end{equation}

We also define

\begin{equation}
S_{1}=X_{1}-X_{2}, \qquad  S_{2}=Y_{1}-Y_{2}
\end{equation}

\begin{equation}
S_{3}=(X_{1}X_{2}Y_{1}+ X_{1}Y_{2}X_{1}+ Y_{1}X_{2}X_{1})- (X_{2}X_{1}Y_{2}+
X_{2}Y_{1}X_{2}+ Y_{2}X_{1}X_{2}) 
\end{equation}

\begin{equation}
S_{4}=(Y_{1}Y_{2}X_{1}+ Y_{1}X_{2}Y_{1}+ X_{1}Y_{2}Y_{1})- (Y_{2}Y_{1}X_{2}+
Y_{2}X_{1}Y_{2}+ X_{2}Y_{1}Y_{2})
\end{equation} 

\begin{equation}
S_{5}=(X_{1}X_{2}X_{1}- X_{2}X_{1}X_{2}), \qquad S_{6}=(Y_{1}Y_{2}Y_{1}-
Y_{2}Y_{1}Y_{2})
\end{equation} 

Using these definitions and the properties $(6.62)$ one obtains quite
generally,

\begin{equation}
\hat{R}_{12}\hat{R}_{23}\hat{R}_{12} -
\hat{R}_{23}\hat{R}_{12}\hat{R}_{23} =
(v+v^{2})S_{1}+(w+w^{2})S_{2}+v^{2}wS_{3}+vw^{2}S_{4}+v^{3}S_{5}+w^{3}S_{6}
\end{equation}

The $S$'s depend on $q$ {\it only}. The dependence on $(v,w)$ are
explicitly displayed in the coefficients of $(6.67)$.

 Now we exploit systematically the constraints on the $S$'s provided by the
known solutions, namely, $(2.32),(2.34),(2,35),(2.36)$ and $(2.37)$.

  The r.h.s. of $(6.67)$ must vanish for the braid solutions $(2.32)$,
namely, for 

\begin{equation}
v=-(1+q^{-2}) \quad w= -(1-q^{-3})
\end{equation}

and also for

\begin{equation}
v=-(1+q^{2}) \quad w= -(1-q^{3})
\end{equation}

Implementing these (for $q \neq 1$) we choose to express $(S_{3},S_{4})$ as

\begin{equation}
S_{3}=\frac {q^{2}}{(q^{2}+1)(q^{2}+q+1)}S_{1}-\frac
{q^{2}}{{(q^{2}+1)}^{2}}S_{2}-\frac
{(q^{2}+1)}{(q^{2}+q+1)}S_{5}-\frac
{{(q^{3}-1)}^{2}}{q{(q^{2}+1)}^{2}}S_{6}
\end{equation}

\begin{equation}
S_{4}=\frac {q^{3}}{{(q^{3}-1)}^{2}}S_{1}+\frac
{q^{2}}{(q^{2}+1)(q^{2}+q+1)}S_{2}-q\frac
{{(q^{2}+1)}^{2}}{{(q^{3}-1)}^{2}}S_{5}+\frac
{(q^{3}-1)(q-1)}{q(q^{2}+1)}S_{6}
\end{equation}

Now note that implementing $(6.62)$ one can express the r.h.s. of $(2.34)$
as 

\begin{equation}
 l_{1}(\hat{R}_{12}- \hat{R}_{23})
+l_{2}({\hat{R}_{12}}^{2}- {\hat{R}_{23}}^{2}) =
(l_{1}v+l_{2}(2v+v^{2}))S_{1}+ (l_{1}w+l_{2}(2w+w^{2}))S_{2}
\end{equation}

Combining this result with $(Case.1)$ or $(2.35)$ one obtains ( since
$v=0$, $l_{2}=0$ )
$$
l_{1}wS_{2}=
\Bigl((w+w^{2})+\frac{q^{2}w^{3}}{{(q^{2}+q+1)}^{2}}\Bigr)S_{2}
$$
\begin{equation}
=\bigl(w+w^{2}\bigr)S_{2}+w^{3}S_{6}
\end{equation}

Hence 
\begin{equation}
S_{6} = q^{2}(q^{2}+q+1)^{-2}S_{2}
\end{equation}

From $(6.70),(6.71)$ and $(6.74)$ one obtains

\begin{equation}
S_{3}=\frac {q^{2}}{(q^{2}+1)(q^{2}+q+1)}S_{1}-\frac
{q(q^{2}-q+1)}{{(q^{2}+1)}^{2}}S_{2}-\frac
{(q^{2}+1)}{(q^{2}+q+1)}S_{5}
\end{equation}

\begin{equation}
S_{4}=\frac {q^{3}}{{(q^{3}-1)}^{2}}S_{1}+\frac
{q(q^{2}-q+1)}{(q^{2}+1)(q^{2}+q+1)}S_{2}-q\frac
{{(q^{2}+1)}^{2}}{{(q^{3}-1)}^{2}}S_{5}
\end{equation}

Combining all the preceding results one finally obtains
 
\begin{equation}
\hat{R}_{12}\hat{R}_{23}\hat{R}_{12} -
\hat{R}_{23}\hat{R}_{12}\hat{R}_{23} =
c_{1}S_{1}+c_{2}S_{2}+c_{5}S_{5}
\end{equation}
$$
 = c_{1}(P^{(-)}_{(12)}-P^{(-)}_{(23)}) +
c_{2}(P^{(0)}_{(12)}-P^{(0)}_{(23)})
$$
$$ 
+c_{5}(P^{(-)}_{(12)}P^{(-)}_{(23)}P^{(-)}_{(12)}-P^{(-)}_{(23)}P^{(-)}_{(12)}P^{(-)}_{(23)})
$$

 Here 

$$ c_{1}=v+v^{2}+\frac {v^{2}wq^{2}}{(q^{2}+1)(q^{2}+q+1)}+vw^{2}{\frac
{q^{3}}{(q^{3}-1)^{2}}}$$
        
 $$ c_{2}=w+w^{2}+\frac
{w^{2}vq(q^{2}-q+1)}{(q^{2}+1)(q^{2}+q+1)}-wv^{2}\frac
{q(q^{2}-q+1)}{(q^{2}+1)^{2}} +\frac
{w^{3}q^{2}}{(q^{2}+q+1)^{2}}$$  
\begin{equation}
c_{5}=\frac{v}{(q^{3}-1)^{2}}\Bigl((q^{3}-1)v+(q^{2}+1)w)\Bigr)
\Bigl((q^{3}-1)v-q(q^{2}+1)w \Bigr)
\end{equation}

As checks one verifies that $c_{1}=c_{2}=c_{5}=0$ for $(6.68)$ and
$(6.69)$. Moreover, $c_{5}=0$ for
$$(q^{2}+1)w=-(q^{3}-1)v$$  
and
$$q(q^{2}+1)w=(q^{3}-1)v$$

Thus one gets back, respectively, $Case 2$ of $(2.36)$ and $Case 3$ of
$(2.37)$.

The form $(6.77)$ makes the dependence on $(v,w)$ entirely explicit, the
$c$'s being given by $(6.78)$ and the $S$'s depending only on $q$. This is
particularly suitable for our purpose. The $P$'s can be reexpressed in
terms of $\hat{R}(v,w)$ using $(1.20)$. But only for $c_{5}=0$ ( cases
$(1,2,3)$ of $(2.35),(2.36),(2.37)$ respectively ) one obtains a relatively
simple form as in $(2.34)$. For $S_{5}$ there is no crucial simplification
as for $S_{6}$ in $(6.74)$. Setting $w=0$ in $(6.77)$ one gets no new
simplification but an identity.

Throughout this paper the projectors $P^{(-)}$ and $P^{(0)}$ have served as
the essential building blocks. They can be easily extracted comparing
$(2.30)$ and $(2.31)$. But for completeness and convenience they are
presented below explicitly.

\begin{equation}
(q^{2}+1)P^{(-)}=
\pmatrix{
0 &0 &0 &0 &0 &0 &0 &0 &0 \cr
0 &1 &0 &-q &0 &0 &0 &0 &0 \cr
0 &0 &q &0 &(q-1)\sqrt{q} &0 &-q &0 &0 \cr
0 &-q &0 &q^{2} &0 &0 &0 &0 &0 \cr
0 &0  &(q-1)\sqrt{q} &0 &(q-1)^{2} &0 &-(q-1))\sqrt{q} &0
 &0 \cr
0 &0 &0 &0 &0  &1 &0 &-q &0 \cr
0 &0 &-q &0 &-(q-1)\sqrt{q} &0 &q &0 &0 \cr
0 &0 &0 &0 &0 &-q &0 &q^{2} &0  \cr
0 &0 &0 &0 &0 &0 &0 &0 &0                                                  
 }.
\end{equation}

\begin{equation}
(q^{2}+q+1)P^{(0)}=
\pmatrix{
0 &0 &0 &0 &0 &0 &0 &0 &0 \cr
0 &0 &0 &0 &0 &0 &0 &0 &0 \cr
0 &0 &1 &0 &\sqrt{q}&0 &q &0 &0 \cr
0 & &0 &0 &0 &0 &0 &0 &0 \cr
0 &0  &\sqrt{q} &0 &q &0 &q\sqrt{q} &0
 &0 \cr
0 &0 &0 &0 &0  &0 &0 &0 &0 \cr
0 &0 &q &0 &q\sqrt{q} &0 &q^{2} &0 &0 \cr
0 &0 &0 &0 &0 &0 &0 &0 &0  \cr
0 &0 &0 &0 &0 &0 &0 &0 &0                                                  
 }.
\end{equation}

\section {APPENDIX B: Nontrivial BE for q=1}

We pointed out in $Sec.2$ that for $a=0$ and $b$ satisfying $(2.38)$ one
obtains two solutions of $BE$ (not modified, with vanishing r.h.s.).
Setting 
$$q=1$$ 
in $(2.38)$ one obtains 
$$b^{2}+3b+1=0$$
or
\begin{equation}
b= {\frac{1}{2}}(-3 \pm {\sqrt{5}})\equiv {-e^{\mp {m}}} 
\end{equation}

( It is amusing to note the relation of $b$, or rather that of $-(b+1)$
with the famous Golden Mean i.e. ${\frac{1}{2}}( \sqrt{5} -1)$. )

Now one obtains

\begin{equation}
\hat {R}(\mp {m})= I-3e^{\mp {m}}\hat {P}^{(0)}
\end{equation}

Where the projector $\hat {P}^{(0)}$ is obtained from $P^{(0)}$ by setting
$q=1$. Thus

 $$
3\hat{P}^{(0)}=
\pmatrix{
0 &0 &0 &0 &0 &0 &0 &0 &0 \cr
0 &0 &0 &0 &0 &0 &0 &0 &0 \cr
0 &0 &1 &0 &1&0 &1 &0 &0 \cr
0 & &0 &0 &0 &0 &0 &0 &0 \cr
0 &0  &1 &0 &1 &0 &1 &0
 &0 \cr
0 &0 &0 &0 &0  &0 &0 &0 &0 \cr
0 &0 &1 &0 &1 &0 &1 &0 &0 \cr
0 &0 &0 &0 &0 &0 &0 &0 &0  \cr
0 &0 &0 &0 &0 &0 &0 &0 &0                                                  
 }.
$$

From $(7.82)$ one easily obtains the nontrivial Hecke condition

\begin{equation}
(\hat {R}(\mp {m})-I)(\hat {R}(\mp {m})+(3e^{\mp {m}} -1)I)=0
\end{equation}

Thus $\hat {R}(\mp {m})$ are {\it not} coboundary matrices. They cannot be
obtained by twisting $I$ since

          $$\Bigl(\hat {R}(\mp {m})\Bigr)^{2} \neq I$$

For 
\begin{equation}
q=1, \quad b=0, \quad a=-1
\end{equation}

again ( from $(2.30),(2.36),(2.37)$ ) one obtains $\hat{R}$ satifying $BE$
given by

\begin{equation}
\hat {R}= I-2\hat {P}^{(-)}
\end{equation}
where $\hat {P}^{(-)}$ is obtained from  $ P^{(-)}$ setting $q=1$ in
$(6.79)$.

 But now, in contrast to $(7.83)$, one has

        $$\hat{R}^{2}=I$$

The $R(\mp m)$ satisfying $YB$ ( Yang-Baxter equation) can be obtained from
$(7.82)$ as 

\begin{equation}
R(\mp {m})=P\hat {R}(\mp {m}) =P-3e^{\mp {m}}\hat {P}^{(0)}
\end{equation}

where the matrix $P$ ( permuting the rows $(2,4),(3,7),(6,8)$ ) leaves $\hat
{P}^{(0)}$ invariant.

In Hieterinta's classification $[8]$ of $4\times 4$ $R$ matrices appear
examples without free parameters. Such a case has been studied $[9]$ in
the context of "exotic bialgebras". Here we have obtained $9\times 9$
examples of such matrices.

\bigskip

\bigskip

I discussed many aspects of this paper with M.A.Sokolov. He also provided,
using a program, the explicit expressions of the coefficients in the
solutions in $Sec.2$. It is a pleasure to thank him . I also thank John
Madore for instructive discussions concerning the noncommutative geometric
aspects.

\smallskip
 
\smallskip

\bibliographystyle{amsplain}

\begin{thebibliography}{}


\bibitem{1}  A. Chakrabarti, RTT relations, a modified braid equation and
noncommutative planes; to be published in J.Math.Phys.;
                                                      (math.QA/0009178)  

\bibitem{2} M.Gerstenhaber and A.Giaquinto, Boundary solutions of the
quantum Yang-Baxter equation and solutions in three dimensions
                                                         (q-alg/9710033)

\bibitem{3} N.Yu.Reshetikhin,L.A.Takhtajan and L.D.Faddeev, Quantization
of Lie groups and Lie algebras; Leningrad Math.J.{\bf 1}(1990)193 

\bibitem{4} A.Klimyk and K.Schmudgen, Quantum groups and their
representations (Spinger,1997)

\bibitem{5} J.Madore, An introduction to noncommutative differntial geometry
and its physical applications (C.U.P.,Second edition,1999)

\bibitem{6} B.L.Cerchiai,G. Fiore and J.Madore, Geometrical tools for
quantum Euclidean spaces, (math.QA/0002007)

\bibitem{7}  E.V.Damaskinsky, P.P.Kulish, M.A.Sokolov,
On the structure of coboundary $R$-matrices for classical series,
Zap.Nauch.Sem. POMI, {\bf 269} (2000), pp. 193-206 (in Russian)


\bibitem{8} J.Hietarinta, Solving the two-dimensional constant quantum
Yang-Baxter equation, J.Math.Phys. {\bf 34}(1993) 1725

\bibitem{9} D.Arnaudon, A.Chakrabarti, V.K.Dobrev and S.G.Mihov, Duality
for exotic bialgebras (math/0101160) 


\end{thebibliography}

\end{document}